\input amstex
\documentstyle{amsppt}
%----------------------------------------------------------------
% Title:     A note on a perfect Euler cuboid
% Author:    Ruslan Sharipov
% Comments:  AmSTeX, 8 pages, amsppt style
% MSC-class: 11D41, 11D72
%----------------------------------------------------------------
%           Replacement for output macro definition
%
\catcode`@=11
\redefine\output@{%
  \def\break{\penalty-\@M}\let\par\endgraf
  \ifodd\pageno\global\hoffset=105pt\else\global\hoffset=8pt\fi  
  \shipout\vbox{%
    \ifplain@
      \let\makeheadline\relax \let\makefootline\relax
    \else
      \iffirstpage@ \global\firstpage@false
        \let\rightheadline\frheadline
        \let\leftheadline\flheadline
      \else
        \ifrunheads@ %\let\makefootline\relax
        \else \let\makeheadline\relax
        \fi
      \fi
    \fi
    \makeheadline \pagebody \makefootline}%
  \advancepageno \ifnum\outputpenalty>-\@MM\else\dosupereject\fi
}
\def\Beta{\mathchar"0\hexnumber@\rmfam 42}
\catcode`\@=\active
%----------------------------------------------------------------
\nopagenumbers
\def\negskp{\hskip -2pt}
\def\Img{\operatorname{Im}}
\def\blue#1{#1}

\catcode`#=11\def\diez{#}\catcode`#=6
\catcode`&=11\catcode`&=4
\catcode`_=11\catcode`_=8
%\catcode`~=11\def\volna{~}\catcode`~=\active
\def\mycite#1{\cite{\blue{#1}}\immediate\special{ps:
     ShrHPSdict begin /ShrBORDERthickness 0 def}}
\def\myciterange#1#2#3#4{\cite{\blue{#2#3#4}}\immediate\special{ps:
     ShrHPSdict begin /ShrBORDERthickness 0 def}}
\def\mytag#1{%
    \tag#1}
\def\mythetag#1{\thetag{\blue{#1}}\immediate\special{ps:
     ShrHPSdict begin /ShrBORDERthickness 0 def}}
\def\myrefno#1{\no#1}
\def\myhref#1#2{\blue{#2}\immediate\special{ps:
     ShrHPSdict begin /ShrBORDERthickness 0 def}}
\def\myEarXivlink{\myhref{http://arXiv.org}{http:/\negskp/arXiv.org}}

\def\mytheorem#1{\csname proclaim\endcsname{Theorem #1}}
\def\mytheoremwithtitle#1#2{\csname proclaim\endcsname{Theorem #1#2}}
\def\mythetheorem#1{\blue{#1}\immediate\special{ps:
     ShrHPSdict begin /ShrBORDERthickness 0 def}}
\def\mylemma#1{\csname proclaim\endcsname{Lemma #1}}
\def\mylemmawithtitle#1#2{\csname proclaim\endcsname{Lemma #1#2}}

\def\mycorollary#1{\csname proclaim\endcsname{Corollary #1}}

\def\mydefinition#1{\definition{Definition #1}}

\def\myconjecture#1{\csname proclaim\endcsname{Conjecture #1}}
\def\myconjecturewithtitle#1#2{\csname proclaim\endcsname{Conjecture #1#2}}

%----------------------------------------------------------------
% Cyrillic fonts definition
%\font\eightcyr=wncyr8
%----------------------------------------------------------------
\pagewidth{360pt}
\pageheight{606pt}
\topmatter
\title
A note on a perfect Euler cuboid.
\endtitle
\author
Ruslan Sharipov
\endauthor
\address Bashkir State University, 32 Zaki Validi street, 450074 Ufa, Russia
\endaddress
\email\myhref{mailto:r-sharipov\@mail.ru}{r-sharipov\@mail.ru}
\endemail
\abstract
    The problem of constructing a perfect Euler cuboid is reduced to a single
Diophantine equation of the degree 12.
\endabstract
\endtopmatter
%\loadbold
%\loadeufb
\TagsOnRight
\document

\head
1. Introduction.
\endhead
     An {\it Euler cuboid\/}, named after Leonhard Euler, is a rectangular 
parallelepiped whose edges and face diagonals all have integer lengths. A 
{\it perfect cuboid\/} is an Euler cuboid whose space diagonal is also of 
an integer length.\par
     In 2005 Lasha Margishvili from the Georgian-American High School in 
Tbilisi won the Mu Alpha Theta Prize for the project entitled "Diophantine 
Rectangular Parallelepiped" (see 
\myhref{http://www.mualphatheta.org/Science_Fair/ScienceFairWinners.aspx}
{http:/\negskp/www.mualphatheta.org/Science\_Fair/...}).
He suggested a proof that a perfect Euler cuboid does not exist. However, 
by now his proof is not accepted by mathematical community. The problem of 
finding a perfect Euler cuboid is still considered as an unsolved problem.
The history of this problem can be found in \mycite{1}. Here are some
appropriate references: \myciterange{2}{2}{--}{35}.
\par
\head
2. Passing to rational numbers.
\endhead
\parshape 16 0pt 360pt 0pt 360pt 180pt 180pt 180pt 180pt 
180pt 180pt 180pt 180pt 180pt 180pt 180pt 180pt 180pt 180pt 
180pt 180pt 180pt 180pt 180pt 180pt 180pt 180pt 180pt 180pt 
180pt 180pt 0pt 360pt 
    Let $A_1B_1C_1D_1A_2B_2C_2D_2$ be a perfect Euler cuboid. Its edges are 
presented by positive integer numbers. \vadjust{\vskip 5pt\hbox to 0pt{\kern 
-20pt \includegraphics{cuboid01.eps}\hss}\vskip -5pt}We write this 
fact as
$$
\align
&\hskip -2em
|A_1B_1|=a,\\
&\hskip -2em
|A_1D_1|=b,
\mytag{2.1}\\
&\hskip -2em
|A_1A_2|=c.
\endalign
$$
Its face diagonals are also presented by positive integers (see Fig\.~2.1):
$$
\align
&\hskip -2em
|A_1D_2|=\alpha,\\
&\hskip -2em
|A_2B_1|=\beta,
\mytag{2.2}\\
&\hskip -2em
|B_2D_2|=\gamma.
\endalign
$$
And finally, the spacial diagonal of this cuboid is presented by a positive 
integer:
$$
\hskip -2em
|A_1C_2|=d.
\mytag{2.3}
$$
From \mythetag{2.1}, \mythetag{2.2}, \mythetag{2.3} one easily derives 
a series of Diophantine equations for the integer numbers $a$, $b$, $c$, 
$\alpha$, $\beta$, $\gamma$, and $d$:
$$
\xalignat 2
&\hskip -2em 
a^2+b^2=\gamma^{\kern 0.4pt 2},
&& b^{\kern 0.4pt 2}+c^{\kern 0.4pt 2}=\alpha^2,\\
\vspace{-1.5ex}
&&&\mytag{2.4}\\
\vspace{-1.5ex}
&\hskip -2em 
c^{\kern 0.4pt 2}+a^2=\beta^{\kern 0.4pt 2},
&& a^2+b^{\kern 0.4pt 2}+c^{\kern 0.4pt 2}=d^{\kern 1pt 2}.
\endxalignat
$$
The main goal of this paper is to reduce the equations \mythetag{2.4}
to a single Diophantine equation for some other integer numbers.\par 
     Relying on the last equation \mythetag{2.4}, we introduce the 
following rational numbers:
$$
\xalignat 3
&\hskip -2em 
x_1=\frac{a}{d},
&&x_2=\frac{b}{d},
&&x_3=\frac{c}{d}.
\mytag{2.5}
\endxalignat
$$
The numbers \mythetag{2.5} are the components of a three-dimensional 
unit vector:
$$
\hskip -2em
(x_1)^2+(x_2)^2+(x_3)^2=1.
\mytag{2.6}
$$ 
From the first three equations \mythetag{2.4} one easily derives the
equations
$$
\align
&\hskip -2em 
(x_1)^2+(x_2)^2=(d_3)^2,\\
&\hskip -2em 
(x_2)^2+(x_3)^2=(d_1)^2,
\mytag{2.7}\\
&\hskip -2em 
(x_3)^2+(x_1)^2=(d_2)^2,
\endalign
$$
where the rational numbers $d_1,\,d_2\,d_3$ are given by the following fractions:
$$
\xalignat 3
&\hskip -2em 
d_1=\frac{\alpha}{d},
&&d_2=\frac{\beta}{d},
&&d_3=\frac{\gamma}{d}.
\mytag{2.8}
\endxalignat
$$
The equations \mythetag{2.6}, \mythetag{2.7}, and \mythetag{2.8} lead to
the following theorem.
\mytheorem{2.1} A perfect Euler cuboid does exist if and only if the
equations \mythetag{2.6} and \mythetag{2.7} are solvable in positive 
rational numbers $x_1,\,x_2,\,x_3$ and $d_1,\,d_2,\,d_3$.
\endproclaim
\demo{Proof} The direct proposition of the theorem~\mythetheorem{2.1}
is immediate from the formulas \mythetag{2.4}, \mythetag{2.5}, and
\mythetag{2.8}. Conversely, assume that $x_1,\,x_2,\,x_3$ and 
$d_1,\,d_2,\,d_3$ are positive rational numbers obeying the equations
\mythetag{2.6} and  \mythetag{2.7}. They a presented by some unique 
irreducible fractions with positive integer numerators and denominators:
$$
\xalignat 6
&x_1=\frac{\nu_1}{\delta_1},
&&x_2=\frac{\nu_2}{\delta_2},
&&x_3=\frac{\nu_3}{\delta_3},
&&d_1=\frac{\nu_4}{\delta_4},
&&d_2=\frac{\nu_5}{\delta_5},
&&d_3=\frac{\nu_6}{\delta_6}.
\endxalignat
$$
Let's denote through $d$ the least common multiple of their denominators,
i\.\,e\. 
$$
d=LCM(\delta_1,\delta_2,\delta_3,\delta_4,\delta_5,\delta_6).
$$
Then the following products are positive integer numbers:
$$
\pagebreak
\xalignat 3
&\hskip -2em 
a=x_1\,d,
&&b=x_2\,d,
&&c=x_3\,d,
\quad\\
\vspace{-1.5ex}
&&&\mytag{2.9}\\
\vspace{-1.5ex}
&\hskip -2em 
\alpha=d_1\,d,
&&\beta=d_2\,d,
&&\gamma=d_3\,d.
\quad
\endxalignat
$$
Applying \mythetag{2.6} and \mythetag{2.7} to \mythetag{2.9}, we derive the
equations \mythetag{2.4} for the integer numbers $a$, $b$, $c$, $\alpha$, 
$\beta$, $\gamma$, and $d$.
\qed\enddemo
\head
3. A rational parametrization. 
\endhead
     Combining \mythetag{2.7} and \mythetag{2.6}, we derive the following
equation for $x_1$ and $d_1$:
$$
\hskip -2em
(x_1)^2+(d_1)^2=1.
\mytag{3.1}
$$
Rational solutions of the equation \mythetag{3.1} are parametrized 
by a rational number $u$:
$$
\xalignat 2
&\hskip -2em 
x_1=\frac{2\,u}{1+u^2},
&&d_1=\frac{1-u^2}{1+u^2}.
\mytag{3.2}
\endxalignat
$$
Since both $x_1$ and $d_1$ in \mythetag{3.1} are positive, the parameter $u$ 
satisfies the inequalities:
$$
\hskip -2em
0<u<1. 
\mytag{3.3}
$$\par
      The second equation in \mythetag{2.7} is $(x_2)^2+(x_3)^2=(d_1)^2$. This 
equation can be written in a form quite similar to the equation \mythetag{3.1}:
$$
\hskip -2em
\Bigl(\frac{x_2}{d_1}\Bigr)^2+\Bigl(\frac{x_3}{d_1}\Bigr)^2=1.
\mytag{3.4}
$$
Rational solutions of the equation \mythetag{3.4} are parametrized 
by a rational number $z$:
$$
\xalignat 2
&\hskip -2em 
\frac{x_2}{d_1}=\frac{2\,z}{1+z^2},
&&\frac{x_3}{d_1}=\frac{1-z^2}{1+z^2}.
\mytag{3.5}
\endxalignat
$$
Combining \mythetag{3.2} and \mythetag{3.5}, we derive the formulas
$$
\hskip -2em
\aligned
&x_2=\frac{2\,z\,(1-u^2)}{(1+u^2)\,(1+z^2)},\\
\vspace{2ex}
&x_3=\frac{(1-u^2)\,(1-z^2)}{(1+u^2)\,(1+z^2)}.
\endaligned
\mytag{3.6}
$$
The parameter $z$ in \mythetag{3.5} and \mythetag{3.6} obeys the inequalities 
similar to \mythetag{3.3}:
$$
\hskip -2em
0<z<1. 
\mytag{3.7}
$$\par
\mytheorem{3.1} The formulas \mythetag{3.2} and \mythetag{3.6} constitute a
rational parametrization of the variables $x_1,\,x_2,\,x_3$ and $d_1$ by
means of two parameters $u$ and $z$ obeying the inequalities \mythetag{3.3} 
and \mythetag{3.7}. The equations
$$
\xalignat 2
&\hskip -2em
(x_1)^2+(x_2)^2+(x_3)^2=1,
&&(x_2)^2+(x_3)^2=(d_1)^2
\mytag{3.8}
\endxalignat
$$
are fulfilled identically due to this parametrization. 
\endproclaim
     The proof of the theorem~\mythetheorem{3.1} is pure calculations.\par
\head
4. An extended parametrization. 
\endhead
     Note that the equations \mythetag{3.8} are two of the four equations
\mythetag{2.6} and \mythetag{2.7} providing a perfect Euler cuboid. 
The other two equations are
$$
\xalignat 2
&\hskip -2em
(x_1)^2+(x_2)^2=(d_3)^2,
&&(x_3)^2+(x_1)^2=(d_2)^2.
\mytag{4.1}
\endxalignat
$$
Let's substitute \mythetag{3.2} and \mythetag{3.6} into the first equation
\mythetag{4.1}. As a result we get 
$$
\hskip -2em
(d_3)^2=\frac{4\,(u^2\,z^2+1)\,(u^2+z^2)}{(1+u^2)^2\,(1+z^2)^2}.
\mytag{4.2}
$$
Similarly, substituting \mythetag{3.2} and \mythetag{3.6} into the second 
equation \mythetag{4.1}, we get
$$
(d_2)^2=\frac{((1+u^2)\,(1+z^2)+2\,z(1-u^2))
\,((1+u^2)\,(1+z^2)-2\,z(1-u^2))}{(1+u^2)^2\,(1+z^2)^2}.
\quad
\mytag{4.3}
$$
Relying on \mythetag{4.2}, we define the following two quantities $\xi$ 
and $a$:
$$
\xalignat 2
&\hskip -2em
\xi=u^2\,z^2+1,
&&a=\frac{d_3\,(1+u^2)\,(1+z^2)}{2\,(u^2\,z^2+1)}.
\mytag{4.4}
\endxalignat
$$
Similarly, relying on \mythetag{4.2}, we define other two quantities 
$\zeta$ and $b$:
$$
\hskip -2em
\aligned 
&\zeta=(1+u^2)\,(1+z^2)+2\,z(1-u^2),\\
\vspace{1ex}
&b=\frac{d_2\,(1+u^2)\,(1+z^2)}{(1+u^2)\,(1+z^2)+2\,z(1-u^2)}.
\endaligned
\mytag{4.5}
$$
The formulas \mythetag{4.4} and \mythetag{4.5} are consistent since the 
denominators of the fractions in them are positive. For $a$ and $b$ from
\mythetag{4.2}, \mythetag{4.3}, \mythetag{4.4}, and \mythetag{4.5}, we
derive. 
$$
\gather
\hskip -2em
a^2=\frac{u^2+z^2}{u^2\,z^2+1},
\mytag{4.6}\\
\vspace{1ex}
\hskip -2em
b^2=\frac{(1+u^2)\,(1+z^2)-2\,z(1-u^2)}{(1+u^2)\,(1+z^2)+2\,z(1-u^2)}.
\mytag{4.7}
\endgather
$$\par
     Since $d_2>0$ and $d_3>0$ (see \mythetag{2.8}, \mythetag{2.2} and
\mythetag{2.3}), the quantities $a$ and $b$ are positive. Therefore,
the formulas \mythetag{4.6} and \mythetag{4.7} define two positive 
functions 
$$
\xalignat 2
&\hskip -2em
a=a(u,z),&&b=b(u,z).
\mytag{4.8}
\endxalignat
$$
The domain of the functions \mythetag{4.8} is outlined by the inequalities
\mythetag{3.3} and \mythetag{3.7}:
$$
\hskip -2em
D_{uz}=\{(u,z)\in\Bbb R^2\!:\,0<u<1\text{\ \ and \ } 0<z<1\}.
\mytag{4.9}
$$
The functions \mythetag{4.8} defined in the domain \mythetag{4.9}
constitute a mapping 
$$
\hskip -2em
f\!:\,D_{uz}\to\Bbb R^2.
\mytag{4.10}
$$
Let's denote through $D_{ab}$ the image of the domain $D_{uz}$ 
under the mapping \mythetag{4.10}:
$$
\hskip -2em
D_{ab}=\Img f=f(D_{uz}).
\mytag{4.11}
$$\par
\parshape 18 0pt 360pt 0pt 360pt 180pt 180pt 180pt 180pt 
180pt 180pt 180pt 180pt 180pt 180pt 180pt 180pt 180pt 180pt 
180pt 180pt 180pt 180pt 180pt 180pt 180pt 180pt 180pt 180pt 
180pt 180pt 180pt 180pt 180pt 180pt 0pt 360pt 
     The domain \mythetag{4.11} is shown in Fig\.~4.1. It is an open
triangle with one curvilinear side. \vadjust{\vskip 5pt\hbox to 
0pt{\kern -5pt \includegraphics{cuboid02.eps}\hss}
\vskip -5pt}The curvilinear side of the triangle $D_{ab}$ is the
graph of the function
$$
b(a)=-1+\frac{2}{a+1}.
$$
Using the formulas \mythetag{4.6} and \mythetag{4.7}, one can prove
that the mapping \mythetag{4.10} sets up a bijective correspondence 
of the points of $D_{uz}$ with the points of $D_{ab}$:
$$
f\!:\,D_{uz}\to D_{ab}.
$$
The inverse mapping 
$$
f^{-1}\!:\,D_{ab}\to D_{uz}
$$
establishing the backward correspondence of the points of $D_{ab}$
with those of $D_{uz}$ is given by two algebraic functions 
$$
\xalignat 2
&\hskip -2em
u=u(a,b),&&z=z(a,b).
\mytag{4.12}
\endxalignat
$$\par
      Let's consider the second formula \mythetag{4.4} and the second 
formula \mythetag{4.5}. We can write these two formulas in the following
way:
$$
\hskip -2em
\gathered
d_2=\frac{(1+u^2)\,(1+z^2)+2\,z(1-u^2)}{(1+u^2)\,(1+z^2)}\,b,\\
\vspace{1ex}
d_3=\frac{2\,(u^2\,z^2+1)}{(1+u^2)\,(1+z^2)}\,a.
\endgathered
\mytag{4.13}
$$
Substituting \mythetag{4.12} into the formulas \mythetag{3.2}, \mythetag{3.6}, 
and \mythetag{4.13}, we can represent $x_1,\,x_2,\,x_3$ and $d_1,\,d_2,
\,d_3$ as functions of two variables $a$ and $b$:
$$
\xalignat 3
&\hskip -2em
x_1=x_1(a,b),&&x_2=x_2(a,b),&&x_3=x_3(a,b),\\
\vspace{-1.5ex}
&&&\mytag{4.14}\\
\vspace{-1.5ex}
&\hskip -2em
d_1=d_1(a,b),&&d_2=d_2(a,b),&&d_3=d_3(a,b).
\endxalignat
$$ 
\mydefinition{4.1} The functions \mythetag{4.14} sharing the common domain $D_{ab}$
constitute a parametrization for the problem of a perfect Euler cuboid. They extend
the rational parametrization given by the functions \mythetag{3.2} and \mythetag{3.6}.
\enddefinition
\mytheorem{4.1} The equations \mythetag{2.6} and \mythetag{2.7} providing a perfect 
Euler cuboid are fulfilled identically by the functions \mythetag{4.14}. 
\endproclaim
     The theorem~\mythetheorem{4.1} is analogous to the theorem~\mythetheorem{3.1}. 
One can see that it is already proved by the above considerations.
\head
5. The characteristic equation. 
\endhead
     Unlike \mythetag{3.2} and \mythetag{3.6}, the functions \mythetag{4.14} 
are not explicit. Below we derive an algorithm for evaluating them. For this
purpose let's return back to the formulas \mythetag{4.4}, \mythetag{4.5}, 
\mythetag{4.6}, and \mythetag{4.7}. From \mythetag{4.4} and \mythetag{4.6} we
derive the equations
$$
\xalignat 2
&\hskip -2em
\xi=u^2\,z^2+1,
&&\xi\,a^2=u^2+z^2.
\mytag{5.1}
\endxalignat
$$ 
Similarly, from \mythetag{4.5} and \mythetag{4.7} we derive the equations
$$
\hskip -2em
\aligned
\zeta&=(1+u^2)\,(1+z^2)+2\,z(1-u^2),\\
\vspace{1ex}
\zeta\,b^2&=(1+u^2)\,(1+z^2)-2\,z(1-u^2).
\endaligned
\mytag{5.2}
$$
Subtracting both equations \mythetag{5.1} from each of the equations \mythetag{5.2},
we get
$$
\hskip -2em
\left\{\aligned 
\zeta-\xi\,(1+a^2)&=2\,z\,(1-u^2),\\
\zeta\,b^2-\xi\,(1+a^2)&=-2\,z\,(1-u^2).
\endaligned\right.
\mytag{5.3}
$$
The equations \mythetag{5.3} constitute a system of two linear algebraic
equations with respect to the variables $\xi$ and $\zeta$. Solving them, we derive
$$
\xalignat 2
&\hskip -2em
\xi=\frac{2\,z\,(1-u^2)\,(1+b^2)}{(1-b^2)\,(1+a^2)},
&&\zeta=\frac{4\,z\,(1-u^2)}{1-b^2}.
\mytag{5.4}
\endxalignat
$$\par
     Now let's substitute $\theta$ for $z^2$ into \mythetag{5.1}. Then the equations
\mythetag{5.1} turn to a system of two linear algebraic equations with respect
to the variables $\xi$ and $\theta$:
$$
\hskip -2em
\left\{\aligned
&\xi-u^2\,\theta=1,\\
&\xi\,a^2-\theta=u^2.
\endaligned\right.
\mytag{5.5}
$$ 
Solving the system of linear equations \mythetag{5.5}, we obtain
$$
\xalignat 2
&\hskip -2em
\xi=\frac{(1-u^2)\,(1+u^2)}{1-a^2\,u^2},
&&\theta=\frac{a^2-u^2}{1-a^2\,u^2}.
\mytag{5.6}
\endxalignat
$$
In \mythetag{5.4} and \mythetag{5.6} we have two expressions for $\xi$. Equating
them we derive the following expression for $z$ expressing it through $a$, $b$ 
and $u$:
$$
\hskip -2em
z=\frac{(1+u^2)\,(1-b^2)\,(1+a^2)}{2\,(1+b^2)\,(1-a^2\,u^2)}.
\mytag{5.7}
$$
Substituting \mythetag{5.7} into the second equation \mythetag{5.4}, we derive: 
$$
\hskip -2em
\zeta=\frac{2\,(1+u^2)\,(1-u^2)\,(1+a^2)}{(1+b^2)\,(1-a^2\,u^2)}.
\mytag{5.8}
$$\par
     Note that the formulas \mythetag{5.6}, \mythetag{5.7}, and \mythetag{5.8}
are similar to each other. They express $z$, $\xi$, $\zeta$, and $\theta$
through $a$, $b$, and $u$. But only two of the three variables $a$, $b$, and $u$
are independent. \pagebreak Due to \mythetag{4.12} the variable $u$ is uniquely 
expressed through $a$ and $b$ within the domain $D_{ab}$ shown in Fig\.~4.1. 
In order to evaluate this expression let's recall that we have the following 
equation:
$$
\hskip -2em
\theta=z^2. 
\mytag{5.9}
$$
Applying \mythetag{5.6} and \mythetag{5.8} to \mythetag{5.9}, we write
\mythetag{5.9} as 
$$
\hskip -2em
\frac{(1+u^2)^2\,(1-b^2)^2\,(1+a^2)^2}{4\,(1+b^2)^2\,(1-a^2\,u^2)^2}
=\frac{a^2-u^2}{1-a^2\,u^2}.
\mytag{5.10}
$$
The denominators of the fractions in the equation \mythetag{5.10} are nonzero 
within the domain $D_{ab}$. For this reason it can be brought to a
polynomial equation:
$$
\gathered
u^4\,a^4\,b^4+(6\,a^4\,u^2\,b^4-2\,u^4\,a^4\,b^2-2\,u^4\,a^2\,b^4)
+(4\,u^2\,b^4\,a^2+\\
+\,4\,a^4\,u^2\,b^2-12\,u^4\,a^2\,b^2+u^4\,a^4+u^4\,b^4+a^4\,b^4)
+(6\,a^4\,u^2+6\,u^2\,b^4-\\
-\,8\,a^2\,b^2\,u^2-2\,u^4\,a^2-2\,u^4\,b^2-2\,a^4\,b^2
-2\,b^4\,a^2)+(u^4+b^4+\\
+\,a^4+4\,a^2\,u^2+4\,b^2\,u^2-12\,b^2\,a^2)+(6\,u^2-2\,a^2-2\,b^2)+1=0.
\endgathered
\quad
\mytag{5.11}
$$
\mytheorem{5.1} The equation \mythetag{5.11} defines the function
$u=u(a,b)$ from \mythetag{4.12} in an implicit form. It is called the
characteristic equation.
\endproclaim
\mytheorem{5.2} A perfect Euler cuboid does exist if and only if 
the characteristic equation \mythetag{5.11} has a rational solution
such that $0<u<1$, while $a$ and $b$ are the coordinates of some
point within the open domain $D_{ab}$ shown in Fig\.~4.1.
\endproclaim
     The inhomogeneous polynomial equation \mythetag{5.11} can
be transformed to a homogeneous equation by adding one more variable $c$:
$$
\gathered
u^4\,a^4\,b^4+6\,a^4\,u^2\,b^4\,c^2-2\,u^4\,a^4\,b^2\,c^2-2\,u^4\,a^2
\,b^4\,c^2+4\,u^2\,b^4\,a^2\,c^4+\\
+\,4\,a^4\,u^2\,b^2\,c^4-12\,u^4\,a^2\,b^2\,c^4+u^4\,a^4\,c^4+u^4\,b^4\,c^4
+a^4\,b^4\,c^4+\\
+\,6\,a^4\,u^2\,c^6+6\,u^2\,b^4\,c^6-8\,a^2\,b^2\,u^2\,c^6-2\,u^4\,a^2\,c^6
-2\,u^4\,b^2\,c^6-\\
-\,2\,a^4\,b^2\,c^6-2\,b^4\,a^2\,c^6+u^4\,c^8+b^4\,c^8+a^4\,c^8
+4\,a^2\,u^2\,c^8+\\
+\,4\,b^2\,u^2\,c^8-12\,b^2\,a^2\,c^8+6\,u^2\,c^{10}-2\,a^2\,c^{10}
-2\,b^2\,c^{10}+c^{12}=0.
\endgathered
\quad
\mytag{5.12}
$$
\mytheorem{5.3} A perfect Euler cuboid does exist if and only if 
the Diophantine equation \mythetag{5.12} has a solution such that $c>0$
and $0<u/c<1$, while $a/c$ and $b/c$ are the coordinates of some
point within the open domain $D_{ab}$ shown in Fig\.~4.1.
\endproclaim
     The theorems~\mythetheorem{5.2} and \mythetheorem{5.3} constitute
the main result of this paper. They can be used in numeric search for a
perfect Euler cuboid.

\enddocument
\end